\documentclass{article}
\usepackage{spconf,amsmath,graphicx}

\usepackage[utf8]{inputenc}
\usepackage[T1]{fontenc}

\usepackage{algorithmic}
\usepackage{algorithm}

\usepackage{hyperref}
\usepackage{amsmath}
\usepackage{amsfonts}
\usepackage{amssymb}
\usepackage{amsthm}
\usepackage{}
\usepackage{graphicx}
\usepackage{color}
\usepackage{xcolor}
\usepackage{diagbox}
\usepackage{url}%\usepackage{slashbox}
\renewcommand{\epsilon}{\varepsilon}
\renewcommand{\phi}{\varphi}
\newcommand*{\N}{\mathbb{N}}

\newcommand*{\R}{\mathbb{R}}

\newcommand*{\C}{\mathbb{C}}

\newcommand{\norm}[1]{\left \| #1 \right \|}
\newcommand{\abs}[1]{\left|#1 \right|}

\theoremstyle{definition}

\theoremstyle{plain}

\usepackage{cite}

\makeatletter
\def\blfootnote{\xdef\@thefnmark{}\@footnotetext}
\makeatother

\begin{document}

% Title.
% ------
\title{FASTER THAN FAST: ACCELERATING THE GRIFFIN-LIM ALGORITHM}
%
% Single address.
% ---------------
\name{Rossen Nenov$^{ \dagger}$ \qquad Dang-Khoa Nguyen$^{\star}$ \qquad Peter Balazs$^{\dagger}$}

\address{$^{\dagger}$ Acoustics Research Institute, Vienna, Austria \\ $^{\star}$ University of Vienna, Austria} 
%\ninept

\maketitle
\begin{abstract}
The phase retrieval problem is found in various areas of applications of engineering and applied physics. 
It is also a very active field of research in mathematics, signal processing and machine learning.
In this paper, we present an accelerated version of the well known Fast Griffin-Lim algorithm (FGLA) for the phase retrieval problem in a general setting. It has increased the speed of convergence, and most importantly, the limit points of the generated sequence can reach a significantly smaller error than the ones generated by FGLA. We will give a motivation of the acceleration and compare it numerically to its predecessors and other algorithms typically used to solve similar problems.
\end{abstract}
\section{Introduction}
\blfootnote{Copyright 2023 IEEE. Published in ICASSP 2023 – 2023 IEEE International Conference on Acoustics, Speech and Signal Processing (ICASSP), scheduled for 4-9 June 2023 in Rhodes Island, Greece. Personal use of this material is permitted. However, permission to reprint/republish this material for advertising or promotional purposes or for creating new collective works for resale or redistribution to servers or lists, or to reuse any copyrighted component of this work in other works, must be obtained from the IEEE. Contact: Manager, Copyrights and Permissions / IEEE Service Center / 445 Hoes Lane / P.O. Box 1331 / Piscataway, NJ 08855-1331, USA. Telephone: + Intl. 908-562-3966.}
\label{sec:intro}

  The \emph{Griffin-Lim algorithm} (GLA) is a famous method to solve the phase retrieval problem in the time-frequency setting \cite{GLA}. It aims to solve two related but different problems, namely constructing a signal from a valid spectrogram with no phase information (Phase recovery) and constructing a signal from a synthetic or modified \emph{Short-Time Fourier transform} (STFT) magnitude. These problems can be seen as an optimization task and many different algorithms from convex { and nonconvex} optimization have been applied {to them} \cite{conv1,conv2,conv3,conv4}. 
  
  Beyond acoustics, phase retrieval problems are found in many different fields involving transformations of signals or data \cite{others1,others2}. In a general setting, the phase retrieval problem can be formulated as the task to find elements in the intersection of the linear subspace $C_1$, i.e. the transformation domain, and of the nonconvex set $C_2$, i.e. those coefficients having a fixed magnitude. The GLA projects the iterates alternatingly onto $C_1$ and onto $C_2$. In \cite{FGLA} applying the \emph{Fast iterative shrinkage-threshold algorithm} (FISTA) was considered for this problem, which resulted in the formulation of FGLA. It outperforms GLA in numerical experiments in speed of convergence and converges to points with smaller error. Since its introduction FGLA was applied for the phase retrieval in many different types of applications \cite{FGLAapp1,FGLAapp2,FGLAapp3}.
  
  Over time different variations of FISTA have been designed for problems beyond the convex setting. This inspired us to use a more flexible algorithm that better fitted for this problem setting and apply it to phase retrieval. It resulted in the formulation of the \emph{Accelerated Griffin-Lim algorithm} (AGLA), which outperforms the established FGLA.
\section{THE PHASE RETRIEVAL PROBLEM}
\label{sec:pagestyle}
Let $T$ be a matrix in $\C^{M\times L}$ with full column rank. We will call $T$ the transformation matrix, like the discrete short time Fourier transform \cite{STFT} or the analysis operator of a finite frame \cite{FrameTh}. Let $s\in [0,\infty)^M$  be the vector of given magnitudes. The phase retrieval problem can be expressed as finding a signal $x^*\in \C^L$, such that the magnitude of the coefficients of the transformation of $x^*$, namely $\abs{Tx^*}$ are as close as possible to $s$, where $\abs{\cdot}$ is understood componentwise. 
{The STFT and other finite frames used in applications are redundant, i.e. $M>L$ and $T$ maps to a true subspace of the coefficient space $\C^M$. Therefore, a nontrivial dependence between the magnitudes and the phases arises, which make the questions about feasibility, uniqueness and stability of the phase retrieval problem difficult to answer. These problems have their own interest and they have been studied in numerous works such as \cite{Grohs, Balan}. 
For results regarding critically sampled transforms, i.e. $M=L$, such as the Discrete Fourier Transform we refer to \cite{oppenheim1975digital}.}

The phase retrieval problem can be formulated  mathematically  in the following form: The signal $x^*$ is the solution of 
      \begin{align}\label{Problem1}
          \min_{x\in \C^L} \norm{\abs{Tx}-s}_2.
      \end{align}
Following the formulation of \cite{FGLA}, this problem can be {rewritten as finding} the intersection or the closest points of two sets.
  
           By $C_1$ we denote the set of possible sequences of coefficients of our transformation, i.e. the range of the transformation matrix $T$:
          \begin{align*}
              C_1=  \{c \in \C^{M} \mid \exists x \in \C^L: c=Tx \}.
          \end{align*}
         This set is a linear subspace and therefore convex. Using the pseudo-inverse $T^\dagger$, which is well defined, the orthogonal projection onto $C_1$ can be expressed as
                  \begin{align*}
              P_{C_1}(c)=TT^\dagger c.
          \end{align*}
          By $C_2$ we denote the set of coefficients, whose magnitude coincide with the given magnitude measurements:
          \begin{align*}
              C_2=  \{c \in \C^{M} \mid \abs{c_i}=s_i \quad \forall i \in \{1,\dots,M\}\}.
          \end{align*}
          This set is compact by definition.
          The projection $P_{C_2}$ onto $C_2$ is equivalent to scaling the elements of $c$ to have magnitude $s$, without changing the phase, which can be written as
          \begin{align*}
              (P_{C_2}(c))_i=\begin{cases} \frac{s_ic_i}{\abs{c_i}} & \text{ if }c_i\neq 0, \\
              s_i & \text{ if }c_i=0. \end{cases}
          \end{align*}
           %One can verify that $\min_{y\in C_{2}} \norm{c-y} = \norm{c-P_{C_2}(c)}$ just as for orthogonal projections holds.
           One can verify that $P_{C_2}(c)$ maps $c$ onto $C_2$ and gives its closest point in $C_2$.
           %One can verify that $P_{C_2}(c)$ is an element in $C_2$ which is nearest to $c$.
          The problem \eqref{Problem1} can be seen as the task of finding coefficients in $C_1$ that are as close as possible to $C_2$, which motivates rewriting the objective as
            \begin{align}\label{problem2}
                \min_{c \in \C^{M}} \delta_{C_1}(c)+ \frac{1}{2}d_{C_2}^2(c),
            \end{align}
            where the indicator function $\delta_C$ of a set $C$ is defined as
          \begin{align*}
              \delta_C(c)=\begin{cases} 0 & \text{if } c \in C, \\
              +\infty & \text{else,}
              \end{cases}
          \end{align*}
           and the distance function $d_C$ to a compact set $C$ is defined as
          \begin{align*}
              d_C(c)=\min_{y\in C } \norm{c-y}.
          \end{align*}
          Solving the problem \eqref{Problem1} is equivalent to finding a solution $c^*\in \C^M$ of \eqref{problem2} since $x^*=T^\dagger c^*$ is then a solution to \eqref{Problem1}.

\section{GRIFFIN-LIM BASED ALGORITHMS}
\label{sec:pagestyle}
{The algorithms will be initialized by a vector $c_0$ {$\in \C^M$ which matches} the magnitudes of $s$ {and the} phases of each component can be {zero or chosen freely}. {We} fix an {integer} $N$, which denotes the amount of iterations of the algorithms.}

The paper \cite{GLA} aims to solve the phase retrieval problem via alternating projections and proposed the following algorithm:

 \begin{figure}[H]
\vspace{-20pt}

\begin{algorithm}[H]

\caption{Griffin-Lim algorithm }\label{alg:GLA}
\begin{algorithmic}
\STATE 
\STATE {\textsc{Fix }} phase $\angle c_0$ and $N \in \N$
\STATE {\textsc{Initialize }}$c_0=s\cdot e^{i \cdot \angle c_0}$ \smallskip
\STATE \hspace{0.5cm}$ \textbf{Iterate for } n=1,\dots, N $
\STATE \hspace{0.5cm}$ c_n =  P_{C_1}(P_{C_2}(c_{n-1})).$
\smallskip 

\STATE {\textsc{Return }} $T^\dagger c_N$
\end{algorithmic}

\end{algorithm}
\vspace{-24pt}
\end{figure}

%\smallskip
%\hrule \vspace{3pt}
%\textbf{Griffin-Lim algorithm (GLA)}
%\hrule \vspace{5pt}
%\noindent 
%\textbf{Fix:} phase $\angle c_0$\\
%\textbf{Initialize:} $c_0=s\cdot e^{i \cdot \angle c_0}$ \\
%\textbf{Iterate:} For $n\geq 1$ 
%\begin{align*}
 %   c_n &=  P_{C_1}(P_{C_2}(c_{n-1})).
%\end{align*}
%\hrule \vspace{4pt}
 The method of alternating projections is a well-known method, which converges in various settings \cite{Altproj},  and for \hyperref[alg:GLA]{GLA} it was proven that the function values of the iterates are nonincreasing and therefore converging \cite{GLA}.
 
 Motivated by the algorithm FISTA \cite{BT}, an inertial step was added in \cite{FGLA}, which resulted in the formulation of \hyperref[alg:FGLA]{FGLA}. It can be seen as \hyperref[alg:GLA]{GLA} but with an enlarged effective stepsize, depending on the last iterates. The case $\alpha=0$ coincides with \hyperref[alg:GLA]{GLA}. The numerical experiments suggest that $\alpha\in [0,1)$ has to hold to have stability in the algorithm, which is coherent with the bounds on the constant $\alpha$ in FISTA for the convex case. 
It seems that the choice $\alpha=0.99$ achieves the best empirical behavior \cite{FGLA}. 

 \begin{figure}[H]
\vspace{-20pt}
\begin{algorithm}[H]
\caption{Fast Griffin-Lim algorithm }\label{alg:FGLA}
\begin{algorithmic}
\STATE 
\STATE {\textsc{Fix }} phase $\angle c_0$, $N \in \N$ and $\alpha>0$
\STATE {\textsc{Initialize }}$c_0=s\cdot e^{i \cdot \angle c_0}$ and $t_0=P_{C_1}(P_{C_2}(c_0))$ \smallskip
\STATE \hspace{0.5cm}$ \textbf{Iterate for } n=1,\dots, N $
\STATE \hspace{0.5cm}$ t_n=P_{C_1}(P_{C_2}(c_{n-1})), $
\STATE \hspace{0.5cm}$ c_n=t_n+\alpha_n (t_n-t_{n-1}). $ \smallskip 

\STATE {\textsc{Return }} $T^\dagger c_N$
\end{algorithmic}

\end{algorithm}
\vspace{-24pt}
\end{figure}

%\smallskip
%\hrule \vspace{3pt}
%\textbf{Fast Griffin-Lim algorithm (FGLA)}
%\hrule \vspace{5pt}
%\noindent 
%\textbf{Fix:} phase $\angle c_0$ and $\alpha>0$\\
%\textbf{Initialize:} $c_0=s\cdot e^{i \cdot \angle c_0}$ and $t_0=P_{C_1}(P_{C_2}(c_0))$ \\
%\textbf{Iterate:} For $n\geq 1$ 
%\begin{align*}
 %   t_n&=P_{C_1}(P_{C_2}(c_{n-1})), \\
  %  c_n&=t_n+\alpha_n (t_n-t_{n-1}). 
%\end{align*} 
%\hrule \vspace{4pt}
  % As we will also see in our numerical results, FGLA appeares to converge faster and to points of less error than GLA. 

Recent developments of new algorithms based on FISTA have lead to many variations of the algorithm, which expanded the applications of inertial methods to certain nonconvex problems \cite{fista1,fista2}. The paper \cite{Laszlo} investigates a relaxation of the FISTA algorithm with two inertial sequences in connection with the minimization of the sum of a convex nonsmooth and a nonconvex differentiable function. %Convergence to critical points was proven for this algorithm under the assumption of the Kurdyka-\L{}ojasiewicz property, a property satisfied by a broad class of functions, i.e. semialgebraic functions, analytic functions and o-minimal functions \cite{kl1,kl2,kl3,kl4}. 

Our function $\delta_{C_1}$ is nonsmooth and convex, while $d_{C_2}^2$ is a nonconvex function and smooth everywhere except on 
\begin{align*}
    \{c\in \C^M \mid \exists i\in \{1,...,{M}\}:  c_i=0 \wedge s_i\neq 0 \}.
\end{align*}
We can still apply the algorithm to solve \eqref{problem2}, but nothing can be deduced from the theoretical statements of previous works. 

\begin{figure}[H]
\vspace{-20pt}
\begin{algorithm}[H]
\caption{Accelerated Griffin-Lim algorithm}\label{alg:AGLA}
\begin{algorithmic}
\STATE 
\STATE {\textsc{Fix }} phase $\angle c_0$, $N\in \N$ and $\alpha,\beta,\gamma>0$
\STATE {\textsc{Initialize }}$c_0=s\cdot e^{i \cdot \angle c_0}$ and $t_0=d_0=P_{C_1}(P_{C_2}(c_0))$ \smallskip
\STATE \hspace{0.5cm}$ \textbf{Iterate for } n=1,\dots, N $
\STATE \hspace{0.5cm}$ t_n = (1-\gamma)d_{n-1}+\gamma P_{C_1}(P_{C_2}(c_{n-1})), $
\STATE \hspace{0.5cm}$ c_n = t_n + \alpha (t_n-t_{n-1}), $ 
\STATE \hspace{0.5cm}$ d_n = t_n + \beta (t_n-t_{n-1}). $\smallskip 

\STATE {\textsc{Return }} $T^\dagger c_N$
\end{algorithmic}

\end{algorithm}
\vspace{-24pt}
\end{figure}

%\smallskip
%%\smallskip \hrule \vspace{3pt}
%textbf{Accelerated Griffin-Lim algorithm (AGLA)}
%\hrule \vspace{5pt}
%\noindent 
%\textbf{Fix:} phase $\angle c_0$ and $\alpha,\beta,\gamma>0$\\
%\textbf{Initialize:} $c_0=s\cdot e^{i \cdot \angle c_0}$ and $t_0=d_0=P_{C_1}(P_{C_2}(c_0))$ \\
%\textbf{Iterate:} For $n\geq 1$ 
%\begin{align*}
 %   t_n &= (1-\gamma)d_{n-1}+\gamma P_{C_1}(P_{C_2}(c_{n-1})), \\
   % c_n &= t_n + \alpha (t_n-t_{n-1}), \\
  %  d_n &= t_n + \beta (t_n-t_{n-1}).
%\end{align*}
%\hrule \vspace{5pt}
%When $\gamma=1$, the proposed AGLA reduces to FGLA. %\tcb{Adding the second inertial sequence and the relaxation parameter $\gamma$ extends the range of $\alpha$ for which the algorithm converges in the numerical experiments, by taking not only the projected of the inertial sequence but also the second inertial sequence into account as well.  Intuitively, this prevents the iterates getting stuck at critical points, if the last step of either the projections or the inertial sequence was sufficiently large.}  Which $\alpha$ one can choose in dependence of $\gamma$ and $\beta$ to theoretically ensure the convergence of the algorithm is beyond the scope of this work. 
{\hyperref[alg:GLA]{GLA} converges to or stops at a point $c^*$ with the property $c^*=P_{C_1}(P_{C_2}(c^*))$ \cite{GLA}, which is not necessarily a global solution. The inertial step of $c_n$ in \hyperref[alg:FGLA]{FGLA} and \hyperref[alg:AGLA]{AGLA} helps surpassing such points, if the last step $t_n-t_{n-1}$ was sufficiently large. Therefore they are expected to converge to points of lower error than \hyperref[alg:GLA]{GLA}. Combining with the nonprojected sequence $d_n$ in \hyperref[alg:AGLA]{AGLA} will enhance this behavior even further, if the distance between the projected $c_n$ and the nonprojected $d_n$ is sufficiently large. {In other words, we are not only taking the direction of the projected iterates into account, but also the direction of the iterates themselves. As can be observed in the numerical experiment, the} additional inertia allows \hyperref[alg:FGLA]{FGLA} and \hyperref[alg:AGLA]{AGLA} to converge to points of lower error, and the relaxation in \hyperref[alg:AGLA]{AGLA} gives more flexibility in the choice of the inertial parameters.}
\section{FURTHER PROJECTION BASED ALGORITHMS}
\label{sec:typestyle}
In this section we present two other algorithms, which also aim to save the phase retrieval problem, to compare them to \hyperref[alg:AGLA]{AGLA}. It was observed in \cite{beyond} that these two in particular perform better than \hyperref[alg:GLA]{GLA} and \hyperref[alg:FGLA]{FGLA}. 

The \emph{Relaxed Averaged Alternating Reflections} (\hyperref[alg:RAAR]{RAAR}) was proposed in \cite{Luke} to solve the phase retrieval problem in Diffraction Imaging. It is a relaxation of the \emph{Douglas-Rachford algorithm}, a well-known method for finding the intersection of two sets using reflections. Especially in the case, where the two sets intersect at a small angle, reflection-based methods converge faster than projection-based methods \cite{beyond}.  Numerical experiments suggest that $\lambda=0.9$ is the best choice for speech signals \cite{beta}. 
\begin{figure}[H]
\vspace{-20pt}
\begin{algorithm}[H]

\caption{Relaxed Averaged Alternating Reflections }\label{alg:RAAR}
\begin{algorithmic}
\STATE 
\STATE {\textsc{Fix }} phase $\angle c_0$, $N\in \N$ and  $0<\lambda \leq 1$
\STATE {\textsc{Initialize }}$c_0=s\cdot e^{i \cdot \angle c_0}$ \smallskip
\STATE \hspace{0.5cm}$ \textbf{Iterate for } n=1,\dots, N $
\STATE \hspace{0.5cm}$ c_{n+1}=\frac{\lambda}{2}\left(c_n+R_{C_1}(R_{C_2}(c_n))\right)+(1-\lambda)P_{C_2}(c_n),$
\STATE \hspace{0.5cm}where $ R_C=2P_{C}-Id $ 
\smallskip 

\STATE {\textsc{Return }} $T^\dagger c_N$
\end{algorithmic}

\end{algorithm}
    \vspace{-24pt}
\end{figure}

%\smallskip 
%\hrule\vspace{5pt}
%\textbf{Relaxed Averaged Alternating Reflections (RAAR)}
%\hrule \vspace{3pt}
%%\noindent 
%\textbf{Fix:} phase $\angle c_0$ and  $0<\lambda \leq 1$\\
%\textbf{Initialize:} $c_0=s\cdot e^{i \cdot \angle c_0}$ \\
%\textbf{Iterate:} For $n\geq 0$ 
%\begin{align*}
  % c_{n+1}=\frac{\lambda}{2}\left(c_n+R_{C_1}(R_{C_2}(c_n))\right)+(1-\lambda)P_{C_2}(c_n),
%\end{align*}
%where $R_C=2P_{C}-Id$. \smallskip 
%\hrule \vspace{5pt}

 We will also compare our algorithm with the \emph{Difference Map} (DM) as proposed in \cite{DM} for phase retrieval problems.
 \begin{figure}[H]
\vspace{-20pt}
\begin{algorithm}[H]
\caption{Difference Map }\label{alg:DM}
\begin{algorithmic}
\STATE 
\STATE {\textsc{Fix }} phase $\angle c_0$, $N\in \N$ and $\rho \in \R\setminus\{0\}$
\STATE {\textsc{Initialize }}$c_0=s\cdot e^{i \cdot \angle c_0}$ \smallskip
\STATE \hspace{0.5cm}$ \textbf{Iterate for } n=1,\dots, N $
\STATE \hspace{0.5cm}$ t_n = P_{C_2}(c_n)+\frac{1}{\rho}(P_{C_2}(c_n)-c_n), $
\STATE \hspace{0.5cm}$ s_n = P_{C_1}(c_n)+\frac{1}{\rho}(P_{C_1}(c_n)-c_n) $ 
\STATE \hspace{0.5cm}$ c_{n+1}=c_n+\rho( P_{C_1}(t_n) - P_{C_2}(s_n)). $\smallskip 

\STATE {\textsc{Return }} $T^\dagger c_N$
\end{algorithmic}

\end{algorithm}
\vspace{-24pt}
\end{figure}
%\hrule\vspace{5pt}
%\textbf{Difference Map (DM)}
%\hrule \vspace{3pt}
%\noindent 
%\textbf{Fix:} phase $\angle c_0$ and $\rho \in \R\setminus\{0\}$\\
%\textbf{Initialize:} $c_0=s\cdot e^{i \cdot \angle c_0}$  \\
%\textbf{Iterate:} For $n\geq 0$ 
%\begin{align*}
%   t_n &= P_{C_2}(c_n)+\frac{1}{\rho}(P_{C_2}(c_n)-c_n), \\
 %  s_n &= P_{C_1}(c_n)+\frac{1}{\rho}(P_{C_1}(c_n)-c_n), \\
  % c_{n+1}&=c_n+\rho( P_{C_1}(t_n) - P_{C_2}(s_n)).
%\end{align*}
%\hrule \vspace{5pt}
The parameter $\rho$ can be chosen arbitrary, but numerical experiments show that the values close to 1 yield the best convergence. In \cite{beyond} the optimal value was observed to be at $\rho=0.8$, which we also observed in our simulations. For the choice $\rho=1$, \hyperref[alg:DM]{DM} is nothing else than \hyperref[alg:RAAR]{RAAR} with $\lambda=1$. It is important to note, that the \hyperref[alg:DM]{DM} algorithm computes twice as many different projections per iteration than the other algorithms presented in this paper. 

\section{NUMERICAL RESULTS}
\label{sec:typestyle}
In this section we will present the results of our numerical experiments and look at the performance of the presented algorithms. As a quality measure we use the \textit{Signal to Noise ratio} (SSNR) of the iterates $c_n$, which is defined as
\begin{align*}
    \text{SSNR}(c)=-10\log_{10}\left(\frac{\norm{\abs{c}-s}}{\norm{s}}\right).
\end{align*}

A reproducible research addendum will be made available at \url{http://bitly.ws/vPr6}, from which one can test different configurations of windows and parameters. 

\begin{figure}[h]
    \centering
    \includegraphics[width=0.47\textwidth]{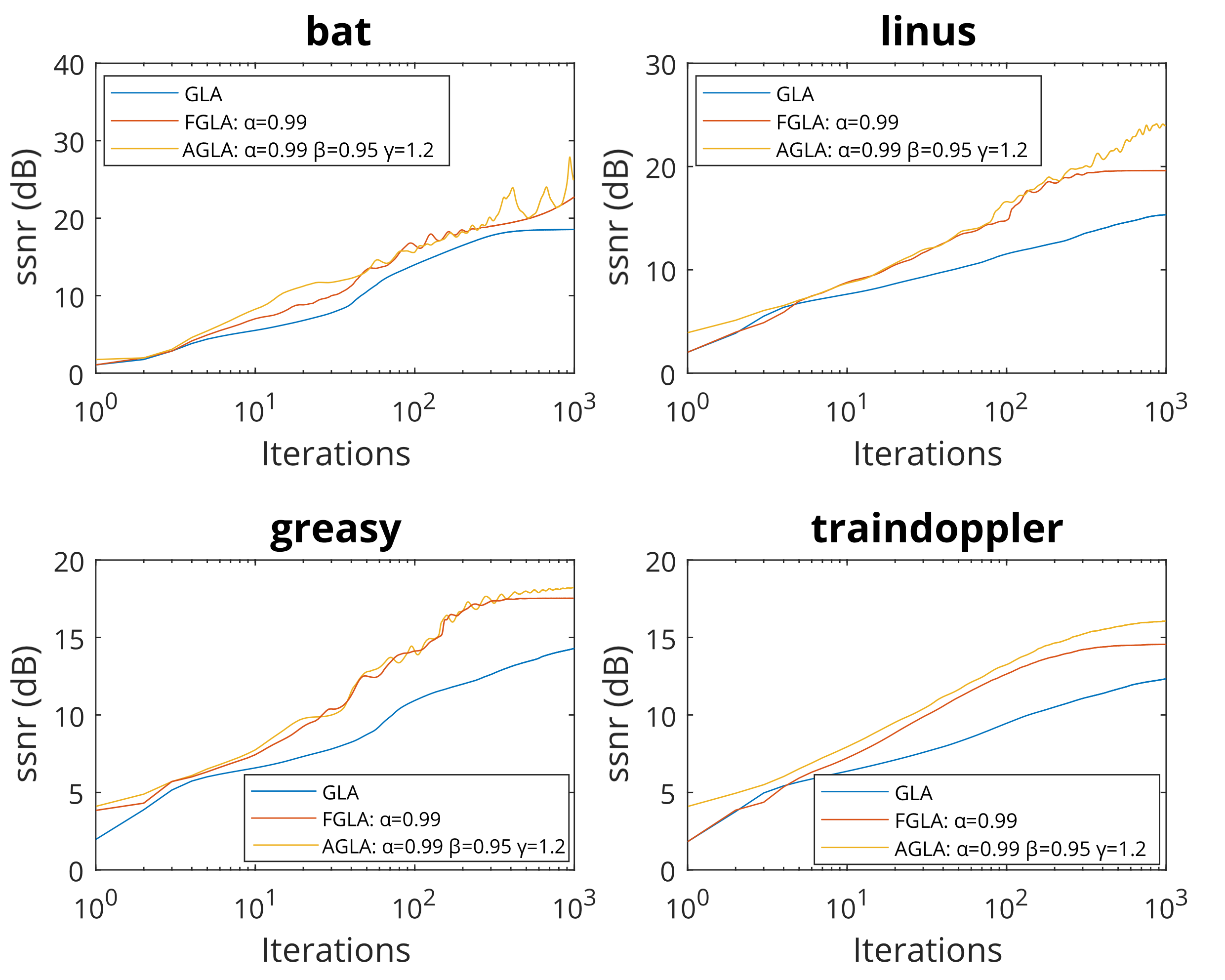}
    \caption{Comparison of the SSNR for the Griffin-Lim based algorithms}
    \label{fig1}
\end{figure} 
\begin{figure}[h!]
    \centering
    \includegraphics[width=0.47\textwidth]{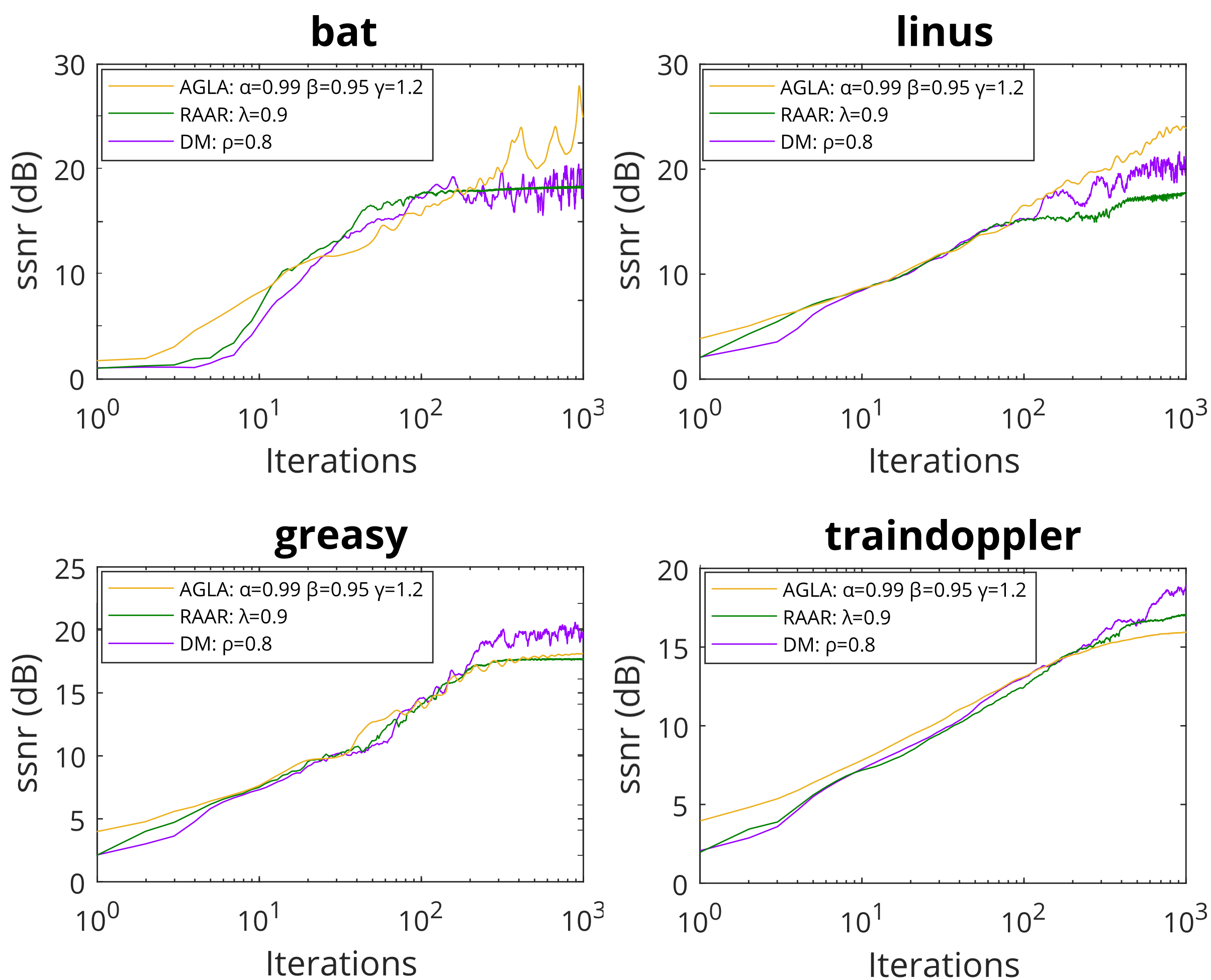}
    \caption{Comparison of the SSNR for \hyperref[alg:AGLA]{AGLA}, \hyperref[alg:DM]{DM} and \hyperref[alg:RAAR]{RAAR}}
    \label{fig2}
\end{figure} 

For the implementation we chose the task of spectogram inversion \cite{Spec} of the STFT. 
The simulations were performed with hop size of 32 and 256 FFT bins using the LTFAT toolbox \cite{LTFAT}. We use a Nuttall window with the window length equaling the number of FFT bins \cite{nuttall}. We computed the spectograms of four different signals, varying in size and noisiness, and applied the algorithms to reconstruct the phase. The initial phase was set to zero and the number of iterations was fixed to $N=1000$. We ran \hyperref[alg:FGLA]{FGLA} with $\alpha=0.99$, \hyperref[alg:RAAR]{RAAR} with $\lambda=0.9$ and \hyperref[alg:DM]{DM} with $\rho=0.8$, based on the previously stated observations of the ideal parameters for convergence. 

\noindent For a fair comparison  we also fixed $\alpha=0.99$, $\beta=0.95$ and $\gamma=1.2$ for \hyperref[alg:AGLA]{AGLA}, which had the best performance on average across our test signals.

In Figure \ref{fig1} we compare \hyperref[alg:GLA]{GLA}, \hyperref[alg:FGLA]{FGLA} and \hyperref[alg:AGLA]{AGLA}. One can notice that \hyperref[alg:AGLA]{AGLA} is an improvement of \hyperref[alg:FGLA]{FGLA} as it converges at least as fast as \hyperref[alg:FGLA]{FGLA} and reaches points with smaller error.

In Figure \ref{fig2} one can see that in some cases \hyperref[alg:AGLA]{AGLA} finds a similar or better critical point than \hyperref[alg:DM]{DM} and \hyperref[alg:RAAR]{RAAR}. For the signal \texttt{traindoppler}, both \hyperref[alg:DM]{DM} and \hyperref[alg:RAAR]{RAAR} pass \hyperref[alg:AGLA]{AGLA} after 300 iterations and converge to a signal with smaller error. %\tcb{For this signal \hyperref[alg:AGLA]{AGLA} improved less compared to \hyperref[alg:FGLA]{FGLA}, which either suggests that depending on the signal some algorithms perform better than the others or that we have not found the best choices for the parameters for \hyperref[alg:AGLA]{AGLA}.}
For this signal and \texttt{greasy} \hyperref[alg:AGLA]{AGLA} improved less compared to \hyperref[alg:FGLA]{FGLA}, which suggests that the optimal choices for the parameter might depend on the nature of the signal.
While \hyperref[alg:DM]{DM} might outperform \hyperref[alg:AGLA]{AGLA} after 300 iterations, we see that the trajectories often times oscillate strongly, whereas \hyperref[alg:AGLA]{AGLA} proceeds smoothly. It is also important to note, that \hyperref[alg:DM]{DM} computes twice as many projections as \hyperref[alg:AGLA]{AGLA} per iteration, which makes it computationally slower, since the projections are the most costly part of these algorithms.

\begin{figure}[t]
\vspace{-10pt}
    \centering
    \includegraphics[width=0.48\textwidth]{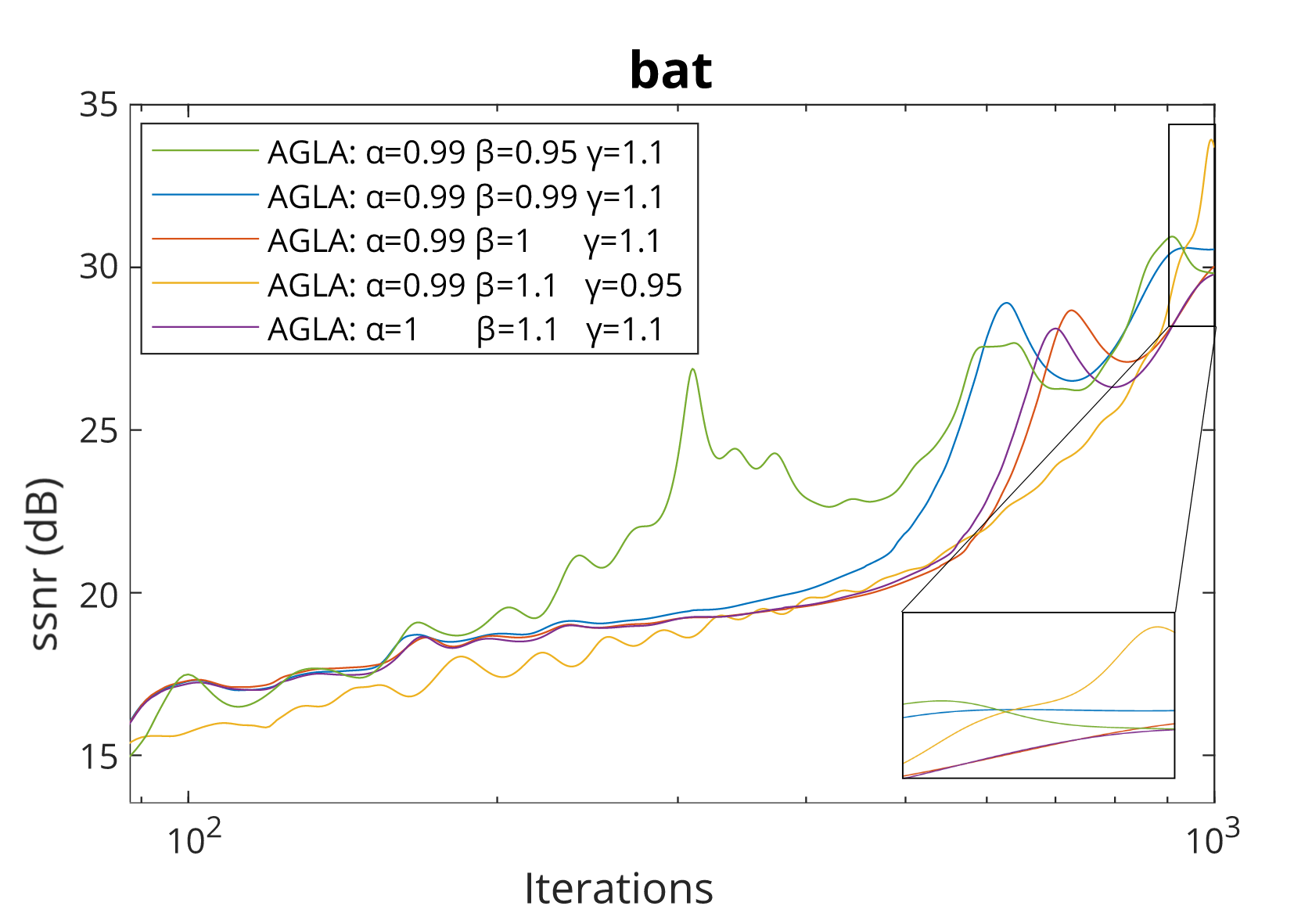}
    \caption{Comparison of the SSNR for the five best parameter choices for \hyperref[alg:AGLA]{AGLA} for the signal \texttt{bat} between 100 and 1000 iterations}
    \label{fig3}
\vspace{-10pt}
\end{figure}

\begin{table}[h]
\centering\scriptsize
\begin{tabular}{|llll|}
\hline
\multicolumn{4}{|c|}{$\alpha=0.99$\rule{0pt}{8pt}} 
\\ \hline    
\multicolumn{1}{|l|}{\backslashbox[18mm]{$\beta$}{$\gamma$}}                       & \multicolumn{1}{c|}{$0.95$}                             & \multicolumn{1}{c|}{$1.1$} & \multicolumn{1}{c|}{$1.2$}  \\ \hline
\multicolumn{1}{|c|}{$0.95$} & \multicolumn{1}{l|}{24.3709}                                                & \multicolumn{1}{l|}{\color{green!40!black} \it 29.8177}                   & 24.9279                   \\ \hline
\multicolumn{1}{|c|}{$0.99$} & \multicolumn{1}{l|}{22.8273}                                                & \multicolumn{1}{l|}{{\color{green!40!black} \textit{30.5563}}}                 & 24.9123                   \\ \hline
\multicolumn{1}{|c|}{$1$}    & \multicolumn{1}{l|}{23.0967}                                                & \multicolumn{1}{l|}{\color{green!40!black} \it 30.0311}                   & 16.4005                   \\ \hline
\multicolumn{1}{|c|}{$1.1$}  & \multicolumn{1}{l|}{\bf {\color{green!40!black} \textit{33.6824}}} & \multicolumn{1}{l|}{21.6791}                   & 16.0700                    \\ \hline \hline
\multicolumn{4}{|c|}{$\alpha =1$\rule{0pt}{8pt}}                                                                                                                                                              \\ \hline
\multicolumn{1}{|l|}{\backslashbox[18mm]{$\beta$}{$\gamma$}}                           & \multicolumn{1}{c|}{$0.95$}                             & \multicolumn{1}{c|}{$1.1$} & \multicolumn{1}{c|}{$1.2$}  \\ \hline
\multicolumn{1}{|c|}{$0.95$} & \multicolumn{1}{l|}{24.8164}                                                & \multicolumn{1}{l|}{20.1219}                   & 17.1421                   \\ \hline
\multicolumn{1}{|c|}{$0.99$} & \multicolumn{1}{l|}{19.6564}                                                & \multicolumn{1}{l|}{20.4847}                   & 15.8881                   \\ \hline
\multicolumn{1}{|c|}{$1$}    & \multicolumn{1}{l|}{21.0041}                                                & \multicolumn{1}{l|}{20.2510}                    & 16.3433                   \\ \hline
\multicolumn{1}{|c|}{$1.1$}  & \multicolumn{1}{l|}{20.2425}                                                & \multicolumn{1}{l|}{\color{green!40!black} \it 29.7887}                   & 16.0721                   \\ \hline\hline
\multicolumn{4}{|c|}{$\alpha =1.1$\rule{0pt}{8pt}}                                                                                                                                                            \\ \hline
\multicolumn{1}{|l|}{\backslashbox[18mm]{$\beta$}{$\gamma$}}                           & \multicolumn{1}{c|}{$0.95$}                             & \multicolumn{1}{c|}{$1.1$} & \multicolumn{1}{c|}{$1.2$}  \\ \hline
\multicolumn{1}{|c|}{$0.95$}  & \multicolumn{1}{l|}{18.1784}                                                & \multicolumn{1}{l|}{13.9603}                   & 13.6818                   \\ \hline
\multicolumn{1}{|c|}{$0.99$}  & \multicolumn{1}{l|}{17.7134}                                                & \multicolumn{1}{l|}{13.2612}                   & 16.8306                   \\ \hline
\multicolumn{1}{|c|}{$1$}     & \multicolumn{1}{l|}{16.9837}                                                & \multicolumn{1}{l|}{14.4684}                   & 14.4627                   \\ \hline
\multicolumn{1}{|c|}{$1.1$}   & \multicolumn{1}{l|}{15.5007}                                                & \multicolumn{1}{l|}{17.6905}                   & 16.7478                   \\ \hline
\end{tabular}
\caption{Comparison of the SSNR for different choices of $\alpha$, $\beta$ and $\gamma$ for the signal \texttt{bat} after 1000 iterations}
\label{tabelle}
\end{table}

In Table \ref{tabelle} we computed the SSNR of the 1000th iterate of the sequence generated by \hyperref[alg:AGLA]{AGLA} with different parameter regimes for one signal {and highlighted the five best choices {for this signal}}. Even though for $\alpha=0.99$, $\beta=1.1$ and $\gamma=0.95$ we achieve the best final error, the overall performance of this regime for the first 500 iterations was significantly worse compared to \hyperref[alg:AGLA]{AGLA} with $\alpha=0.99$, $\beta=0.99$ and $\gamma=1.1$. This phenomenon can be observed in Figure \ref{fig3}. 

For \hyperref[alg:FGLA]{FGLA} it was observed in \cite{FGLA}, that for $\alpha\geq 1$ the algorithm either oscillates strongly, converges to point of great error or completely breaks. \hyperref[alg:AGLA]{AGLA} still performed overall best for $\alpha \in (0,1)$, but interestingly there are cases with $\alpha=1$, where \hyperref[alg:AGLA]{AGLA} results in a small error as well.{
With three parameters in \hyperref[alg:AGLA]{AGLA} it is harder to find the optimal choice of parameters, but they {also} make \hyperref[alg:AGLA]{AGLA} more flexible and allow it to perform well for a wider range of parameters.}

%Overall we see that the addition of the relaxation and the second inertial sequence allow  , but as we see in Figure \ref{fig3}, finding the optimal choice of parameters is more challenging.}

% Please add the following required packag es to your document preamble:

%   In Figure 1 we applied the algorithms to reconstruct a signal from its spectogramm. The initial phase for all algorithm was set to zero and the amount of iterations was $10^3$. For FGLA we chose $\alpha=3$, for AFLGA $\gamma=,\alpha=$ and $\beta=$ and for the relaxed Douglas-Rachford $\lambda=$. What we notice is that AFLGA converges faster than the other algorithms, if you look at the average slope, and has a smaller error. 
% Compared to GLA the other algorithms are not strictly increasing, which is typical for inertial type algorithms. 

\section{CONCLUSION AND FUTURE WORK}
\label{sec:copyright}
In this paper we presented a new algorithm to solve the phase retrieval problem.
Compared to \hyperref[alg:GLA]{GLA} and \hyperref[alg:FGLA]{FGLA}, \hyperref[alg:AGLA]{AGLA} converges in numerical examples faster and to better points with hardly any extra computational cost. We also saw that depending on the signal, \hyperref[alg:AGLA]{AGLA} can perform better than \hyperref[alg:DM]{DM} and \hyperref[alg:RAAR]{RAAR}. While these two algorithms sometimes find better critical points, they tend to have strong oscillations in their function values after a certain amount of iterations compared to \hyperref[alg:AGLA]{AGLA}. This sparks the interest in investigating the performance of hybrid regimes as proposed in \cite{beyond}, where one runs either \hyperref[alg:DM]{DM} or \hyperref[alg:RAAR]{RAAR} for a set amount of iterations to escape the neighbourhood of bad local minima and then applies \hyperref[alg:AGLA]{AGLA} for the rest of the iterations. {Furthermore, the question arises how well \hyperref[alg:AGLA]{AGLA} performs with different initialization, {other} window functions and lower redundancies, and if the optimal choices for the parameters might vary. }

In our future research we will look into stating convergence results of \hyperref[alg:FGLA]{FGLA} and \hyperref[alg:AGLA]{AGLA}, namely proving convergence of the iterates $(c_n)_{n\in\mathbb{N}}$ and of the objective function values $\left( d_{C_2}^2(c_n) \right) _{n\in\mathbb{N}}$. 
%\tcm{As suggested by \cite{Laszlo}, this can be done with the help of the Kurdyka-Łojasiewicz property.} 
% The study of the convergence of these algorithms may also help in finding optimal choices for the parameters $\alpha$, $\beta$ and $\gamma$.

%\vfill
\section{Acknowledgment}
\label{sec:refs}
This work was supported by the Austrian Science Fund FWF-project NoMASP (“Nonsmooth nonconvex optimization methods for acoustic signal processing”; P 34922-N). The authors would like to thank Prof. Radu Ioan Bo\c{t} (University of Vienna) for his encouragement and enthusiastic discussions that initiated this work {and Dr. Nicki Holighaus (Austrian Academy of Sciences) for his valuable input during the reviewing process.}
%\vfill
%\pagebreak
%\vfill \pagebreak \vfill \pagebreak

\bibliographystyle{IEEEtran}

\end{document}